\documentclass[11pt]{article}
\usepackage{mathrsfs}
\usepackage[centertags]{amsmath}
\usepackage{amsfonts}
\usepackage{amssymb}
\usepackage{amsthm}
\usepackage{cases}
\usepackage{indentfirst}
\usepackage{epsfig}
\usepackage{graphicx}
\usepackage{color}
\usepackage{CJK}
\usepackage[numbers,sort&compress]{natbib}
\usepackage{subfigure} 

\date{}

\newtheorem{theorem}{Theorem}[section]

\newtheorem{lemma}{Lemma}[section]

\newtheorem{proposition}{Proposition}[section]
\newtheorem{remark}{Remark}[section]

\numberwithin{equation}{section}

\allowdisplaybreaks[4]

\begin{document}
\title{\textbf{Expansions of maximum and minimum from Generalized Maxwell distribution}
\thanks{
{E-mail address: hjw1303987297@126.com (J. Huang)}
}
\author{Jianwen Huang\\
{\small School of Mathematics $\&$ Statistics, Tianshui Normal University, Tianshui, 741001 China}
}
}

\maketitle
\begin{quote}
{\bf Abstract.}~~Generalized Maxwell distribution is an extension of the classic Maxwell distribution. In this paper, we concentrate on the joint distributional asymptotics of normalized maxima and minima. Under optimal normalizing constants, asymptotic
expansions of joint distribution and density for normalized partial maxima and minima are established. These expansions are used to educe speeds of convergence of joint distribution and density of normalized maxima and minima tending to its corresponding ultimate limits. Numerical analysis are provided to support
our results.

\end{quote}

\section{Introduction}
\label{sec1}

Let independent and identically distributed (i.i.d.) random variables $\{X_i,~i\geq 1\}$ with marginal generalized Maxwell distribution $F_k$ (abbreviated as $F_k\sim \mbox{GMD}(k)$), where $k>0$ is the shape parameter. The probability density function (pdf) of $\mbox{GMD}(k)$ is given by \cite{Voda 2009}
\begin{align}\label{eq.1}
f_k(x)=\frac{k}{2^{k/2}\sigma^{2+1/k}\Gamma(1+k/2)}x^{2k}\exp\left(-\frac{x^{2k}}{2\sigma^2}\right),~x\geq0,
\end{align}
where $\sigma$ is a positive parameter and $\Gamma(\cdot)$ denotes the Gamma function. For $k=1$, $\mbox{GMD}(k)$ degenerates to the classical Maxwell distribution.

In the last few years, asymptotic behaviors regarding the normalized maxima of $\mbox{GMD}(k)$ have been discussed in the literature. Huang and Chen \cite{Huang and Chen 2016} established the Mills' type inequality and Mills' type ratio of $\mbox{GMD}(k)$, and demonstrated that
\begin{align}\label{eq.2}
\lim_{n\to\infty}\mathbb{P}(M_n\leq u(x,b_n))=\Lambda(x)
\end{align}
with
\begin{align}\label{eq.3}
u(x,b_n)=k^{-1}\sigma^2b_n^{1-2k}x+b_n
\end{align}
and the normalizing constants $b_n$ is determined by
\begin{align}\label{eq.4}
1-F_k(b_n)=\frac{1}{n},
\end{align}
where $M_n=\max\{X_i,~1\leq i\leq n\}$ represents the partial maximum of $\{X_i,~i\geq 1\}$ and $\Lambda(x)=\exp(-e^{-x})$ stands for the Gumbel distribution. Additionally, Huang and Wang \cite{Huang and Wang 2018} gained the asymptotic expansions for distribution as follows:
\begin{align}\label{eq.28}
b^{2k}_n\left[b^{2k}_n(\mathbb{P}(M_n\leq u(x,b_n))-\Lambda(x))-l_k(x)\Lambda(x)\right]\to
\left[w_k(x)+\frac{1}{2}l^2_k(x)\right]\Lambda(x)
\end{align}
as $n\to\infty$, where $l_k(x)$ and $w_k(x)$ are respectively given by
\begin{align}\label{eq.5}
l_k(x)=2^{-1}k^{-1}\sigma^2[(2k-1)x^2-2x]e^{-x}
\end{align}
and
\begin{align}\label{eq.6}
w_k(x)=-24^{-1}k^{-2}\sigma^4[3(2k-1)^2x^4-4(2k+1)(2k-1)x^3+24x^2-48kx]e^{-x}.
\end{align}

Related work concerning properties of $\mbox{GMD}(k)$ can be found in \cite{Nadarajah et al 2013} \cite{Liu and Liu 2013} \cite{Dar et al 2017} and \cite{Huang et al 2017 Stat}. $\mbox{GMD}(k)$ has received applications in a variety of areas, see, e.g., establishing the discretized equilibrium distributions of the lattice Boltzmann method \cite{Shim 2017}, describing the stress relaxation behavior of food matrices \cite{Del Nobile et al 2007}, being used to derive a different class of multiple relaxation-time lattice Boltzmann models \cite{Asinari and Karlin 2009}, and modeling the linear viscoelastic fluids \cite{Liu Q et al 2011}.

For the basic properties, which determine the max-domain of attractions the distribution attributes, see \cite{Liao et al 2014a} \cite{Lin et al 2016} \cite{Du and Chen 2016} and \cite{Liao et al 2014b}. In extreme value theory, asymptotics for the joint distribution of maxima and minima is an interesting topic in recent literature. For sequence of univariate random variables under conventional assumptions, see \cite{Leadbetter et al 2012}, \cite{Hashorva and Weng 2014} and \cite{Liao et al 2019}. For bivariate Gaussian arrays under relaxation assumptions, see \cite{Hashorva and Weng 2013}. For dissimilar improved assumptions, see \cite{Liao and Peng 2015} and \cite{Lu and Peng 2017}. Let $m_n=\min\{X_i,~1\leq i\leq n\}$ denote the partial minimum of $\{X_i,~i\geq 1\}$. The target of this paper is to establish asymptotic expansions for distribution
and density of $(M_n,~m_n)$ for $\mbox{GMD}(k)$, from which we can respectively educe the speeds of convergence of the distribution and density from normalized $(M_n,~m_n)$ to the distribution and density of the associating extreme value distribution.

The contents of this paper are constructed as follows. Section \ref{sec2} provides the main results, and numerical analysis is given in Section \ref{sec3}. A few necessary auxiliary lemmas and the proofs of the main results are presented in Section \ref{sec.4}.

\section{Main results}
\label{sec2}

In this section, we present the main results. Proposition \ref{pro.1} provides the asymptotic expansions of a distribution for the minimum of the $\mbox{GMD}(k)$. Theorem \ref{the.1} displays the asymptotic distribution of the joint distribution of normalized $(M_n,~m_n)$. The asymptotic expansion for its distribution is showed in Theorem \ref{the.2}, and the associating asymptotic density expansion is revealed in Theorem \ref{the.3}.

\begin{proposition}
\label{pro.1}
Let $\{X_i,~i\geq 1\}$ be an i.i.d. random sequence with common distribution function $F_k$ obeying the $\mbox{GMD}(k)$. Set $m_n=\min\{X_i,~1\leq i\leq n\}$. Then,
\begin{align}\label{eq.7}
\notag &b^{2k}_n\left\{b^{2k}_n\left[\mathbb{P}(m_n\leq v(y,b_n))-(1-\Lambda(-y))\right]+l_k(-y)\Lambda(-y)\right\}\\
&\to -\left[w_k(-y)+\frac{1}{2}l^2_k(-y)\right]\Lambda(-y),
\end{align}
as $n\to\infty$, where
\begin{align}\label{eq.8}
v(y,b_n)=k^{-1}\sigma^2b^{1-2k}_ny-b_n,
\end{align}
$l_k(\cdot)$, $w_k(\cdot)$ and $b_n$ are respectively given by (\ref{eq.5}), (\ref{eq.6}) and (\ref{eq.4}).
\end{proposition}

\begin{remark}
Noting that $b^{-2k}_n=O((\log n)^{-1})$ in view of (\ref{eq.4}), it follows from Proposition \ref{pro.1} that the convergence rate of the distribution of normalized minima under the normalized constant $b_n$ converging to its extreme value distribution is the same order as $O((\log n)^{-1})$.
\end{remark}

In the following result, we establish the limiting distribution of the joint distribution of normalized $(M_n,~m_n)$ from i.i.d. random variables.

\begin{theorem}\label{the.1}
Let $\{X_i,~i\geq 1\}$ denote i.i.d. random variables with common distribution function $\mbox{GMD}(k)$ $F_k$. Set $M_n=\bigvee^n_{i=1}X_i$ and $m_n=\bigwedge^n_{i=1}X_i$ respectively indicate the maxima and the minima. We have
\begin{align}\label{eq.9}
\mathbb{P}(M_n\leq u(x,b_n),~m_n\leq v(y,b_n))\to\Lambda(x)(1-\Lambda(-y)),
\end{align}
as $n\to\infty$, where $u(x,b_n)$ and $v(y,b_n)$ are respectively determined by (\ref{eq.3}) and (\ref{eq.8}) with the normalizing constant $b_n$ defined by (\ref{eq.4}).
\end{theorem}

\begin{remark}
A combination of (\ref{eq.2}) and (\ref{eq.7}), Theorem \ref{the.1} demonstrates that $M_n$ and $m_n$ are asymptotically independent.
\end{remark}

In the following result, we present asymptotic expansion of a distribution for $(M_n,~m_n)$, by which one derives the corresponding speed of convergence of the joint distribution to its limit. Denote $l_k(x,y)=l_k(x)+l_k(-y)$ and $w_k(x,y)=w_k(x)+w_k(-y)$.

\begin{theorem}\label{the.2}
With $u(x,b_n)$ and $v(y,b_n)$ defined by (\ref{eq.3}) and (\ref{eq.8}), we get
\begin{align}\label{eq.10}
\notag&b^{2k}_n\bigg\{b^{2k}_n\left[\mathbb{P}(M_n\leq u(x,b_n),~m_n\leq v(y,b_n))-\Lambda(x)(1-\Lambda(-y))\right]\\
\notag&-l_k(x)\Lambda(x)+l_k(x,y)\Lambda(x)\Lambda(-y)\bigg\}\\
&\to\left[w_k(x)+\frac{1}{2}l_k^2(x)\right]\Lambda(x)
-\left[w_k(x,y)+\frac{1}{2}l_k^2(x,y)\right]\Lambda(x)\Lambda(-y)
\end{align}
as $n\to\infty$, where $l_k(\cdot)$ and $w_k(\cdot)$ are respectively provided by (\ref{eq.5}) and (\ref{eq.6}).
\end{theorem}

\begin{remark}
By (\ref{eq.4}), one can easily check that $b^{2k}_n=O(\log n)$. Accordingly, Theorem \ref{the.2} clears that the speed of convergence of $\mathbb{P}(M_n\leq u(x,b_n),~m_n\leq v(y,b_n))$ tending to its ultimate limit is proportional to $(\log n)^{-1}$.
\end{remark}

Theorem \ref{the.2} could be used to deduce asymptotic expansions for the density function of $(M_n,~m_n)$. Let
\begin{align}\label{eq.27}
g_n(x,y)=k^{-2}\sigma^4b^{2-4k}_nn(n-1)\left[F_k(u(x,b_n))-F_k(v(y,b_n))\right]^{n-2}f_k(u(x,b_n))f_k(v(y,b_n))
\end{align}
stand for the pdf of the normalized $(M_n,~m_n)$, and let
\begin{align}
\notag \Delta_n(g_n,g;x,y)=g_n(x,y)-g(x,y),
\end{align}
where $g(x,y)=\Lambda(x)\Lambda(-y)e^{-x}e^y$ is the joint density function of $\Lambda(x)\Lambda(-y)$. By means of Proposition $2.5$ in \cite{Resnick 1987}, it results in $\lim_{n\to\infty}\Delta_n(g_n,g;x,y)=0$.

In the following result, we establish the higher
order expansion of the probability density for normalized $(M_n,~m_n)$, from which we obtain the relevant convergence speed of the joint density function to its limiting density.
\begin{theorem}\label{the.3}
For $u(x,b_n)$ and $v(y,b_n)$ determined by (\ref{eq.3}) and (\ref{eq.8}), we get
\begin{align}
\notag&b^{2k}_n\left\{b^{2k}_n\Delta_n(g_n,g;x,y)-C_1(x,y)\Lambda(x)\Lambda(-y)e^{-x}e^y\right\}\to C_2(x,y)\Lambda(x)\Lambda(-y)e^{-x}e^y,
\end{align}
as $n\to\infty$, where
\begin{align}
\notag C_1(x,y)=l_k(x,y)-\sigma^2\left[\frac{2k-1}{2k}(x^2+y^2)-2(x-y)+\frac{2}{k}\right]
\end{align}
and
\begin{align}
\notag C_2(x,y)=&\sigma^4\bigg\{\frac{1}{k}\bigg[\frac{(2k-1)^2}{8k}(x^4+y^4)-\frac{(2k+1)(4k-1)}{3k}(x^3-y^3)
+\frac{(2k+1)(2k-1)}{2k}(x^2+y^2)\\
\notag&-2(x-y)+4\bigg]+\left(\frac{2k-1}{2k}x^2-2x+\frac{1}{k}\right)\left(\frac{2k-1}{2k}y^2+2y+\frac{1}{k}\right)
\bigg\}\\
\notag&-\sigma^2\left[\frac{2k-1}{2k}(x^2+y^2)-2(x-y)+\frac{2}{k}\right]l_k(x,y)+w_k(x,y)
+\frac{1}{2}l^2_k(x,y).
\end{align}
\end{theorem}

\begin{remark}
Observing that $b^{2k}_n=O(\log n)$, Theorem \ref{the.3} evidences that the speed of convergence of the joint pdf of the normalized $(M_n,~m_n)$ tending to its ultimate limit is the same order as $O((\log n)^{-1})$.
\end{remark}

\section{Numerical analysis}
\label{sec3}

In this section, numerical research is presented to demonstrate the precision of higher order expansions for the joint distribution and the joint density of normalized $(M_n,~m_n)$. It follows from (\ref{eq.20}) and (\ref{eq.27}) that the real values of the joint distribution of $(M_n,~m_n)$ are equal to $\mathbb{P}(M_n\leq u(x,b_n),~m_n\leq v(y,b_n))=F^n_k(u(x,b_n))-[F_k(u(x,b_n))-F_k(v(y,b_n))]^n$ and the real values of the joint density of $(M_n,~m_n)$ are provided by (\ref{eq.27}). Let $S_i(x,y)$ and $T_i(x,y)$, $i=1,2,3$, respectively represent the first-order, the second-order and the third-order asymptotics of the joint distribution and the joint density of normalized $(M_n,~m_n)$. Due to Theorems \ref{the.2} and \ref{the.3}, we can obtain
\begin{align}
\notag &S_1(x,y)=\Lambda(x)-\Lambda(x)\Lambda(-y),\\
\notag &S_2(x,y)=S_1(x,y)+b^{-2k}_n[l_k(x)\Lambda(x)-l_k(x,y)\Lambda(x)\Lambda(-y)],\\
\notag &S_3(x,y)=S_2(x,y)+b^{-4k}_n\left\{\left[w_k(x)+\frac{1}{2}l_k^2(x)\right]\Lambda(x)
-\left[w_k(x,y)+\frac{1}{2}l_k^2(x,y)\right]\Lambda(x)\Lambda(-y)\right\},
\end{align}
and
\begin{align}
\notag &T_1(x,y)=\Lambda(x)\Lambda(-y)e^{-x}e^y,\\
\notag &T_2(x,y)=T_1(x,y)+C_1(x,y)\Lambda(x)\Lambda(-y)e^{-x}e^yb^{-2k}_n,\\
\notag &T_3(x,y)=T_2(x,y)+C_2(x,y)\Lambda(x)\Lambda(-y)e^{-x}e^yb^{-4k}_n,
\end{align}
where $b_n$ is defined in (\ref{eq.4}).

\begin{figure}[htbp]
\centering
\subfigure[$(x,y)=(5,5),~\sigma=1$.]{
\includegraphics[scale=0.5]{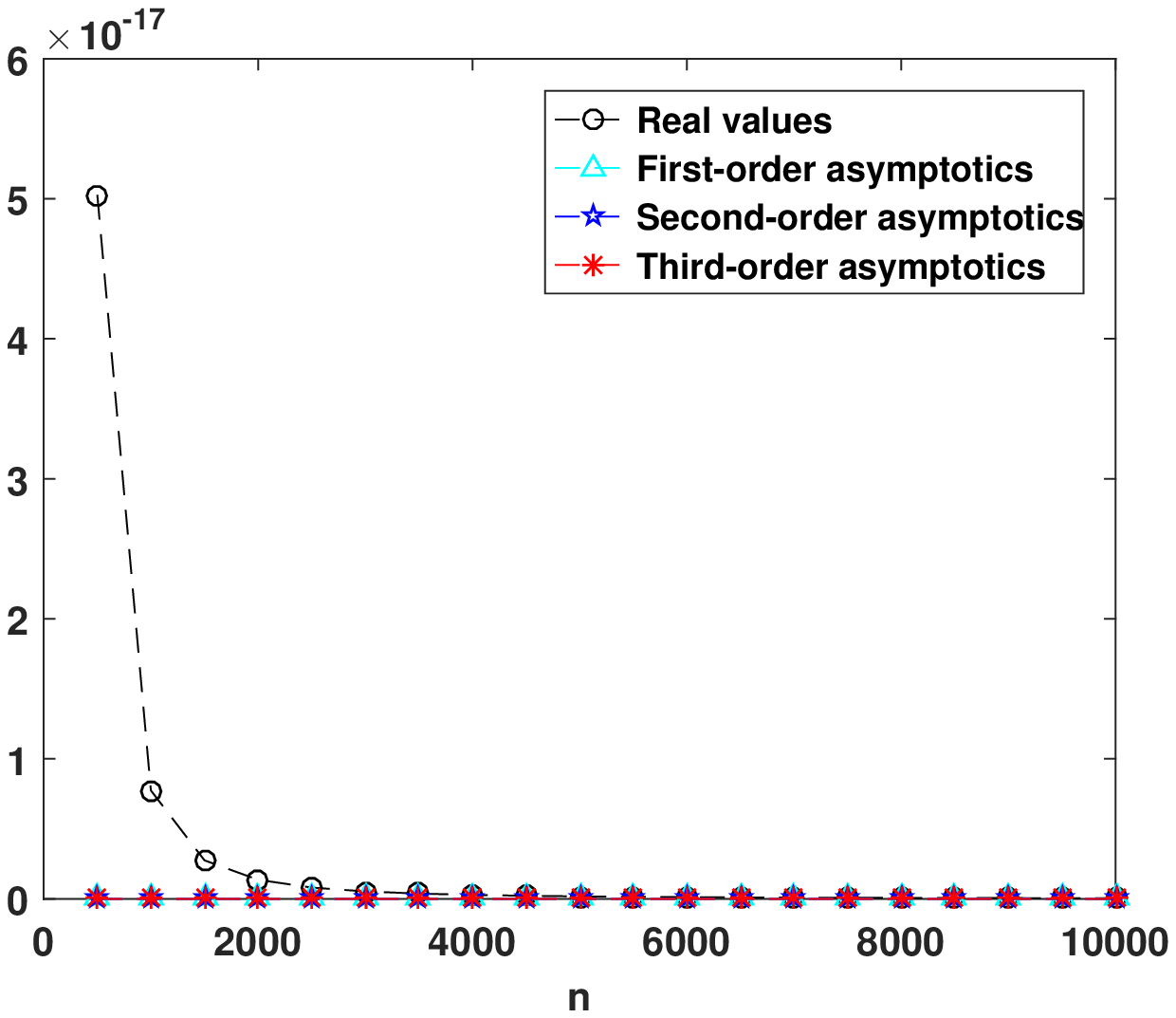}
}
\quad
\subfigure[$(x,y)=(2,6),~\sigma=1$.]{
\includegraphics[scale=0.5]{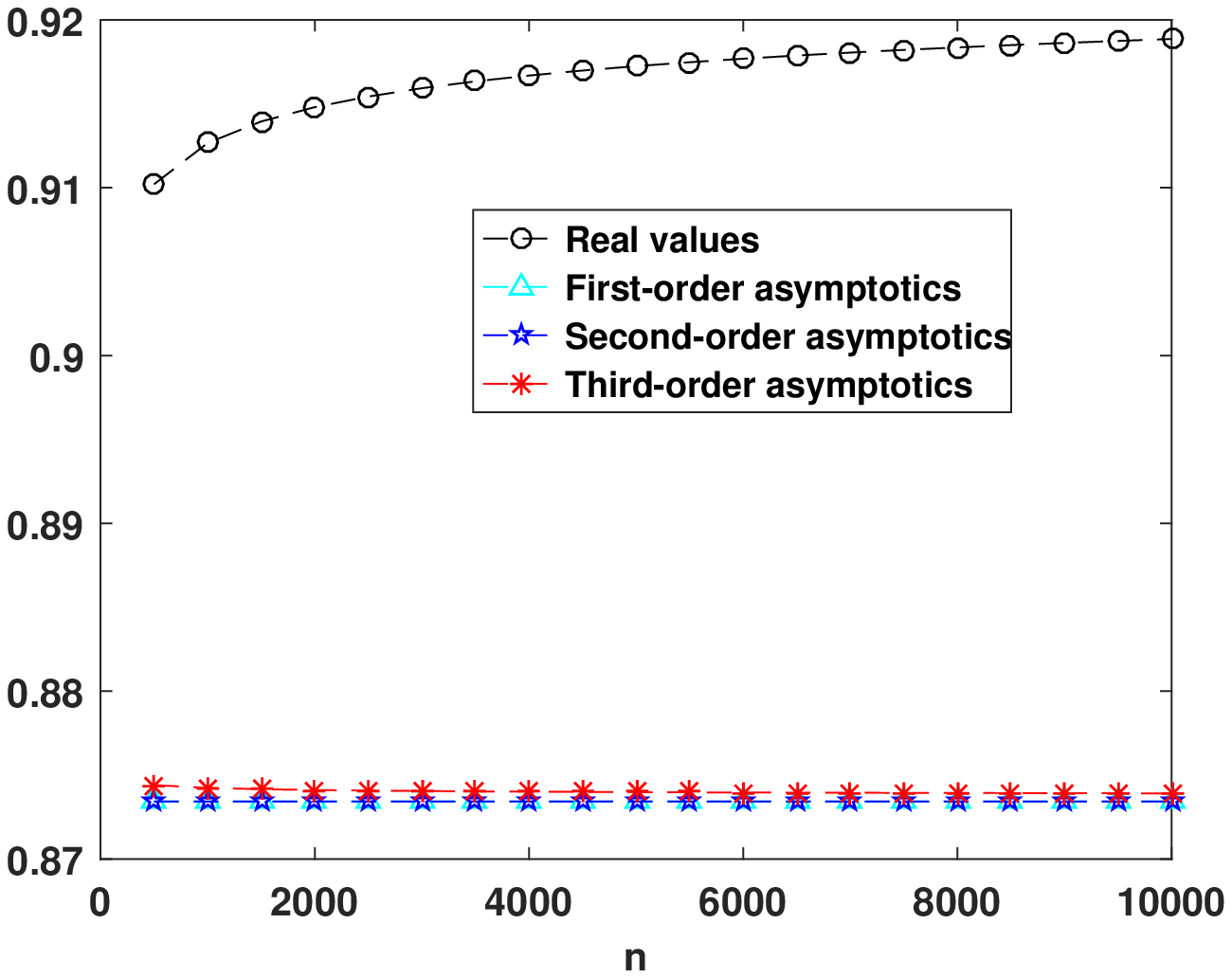}
}
\quad
\subfigure[$(x,y)=(6,2),~\sigma=1$.]{
\includegraphics[scale=0.5]{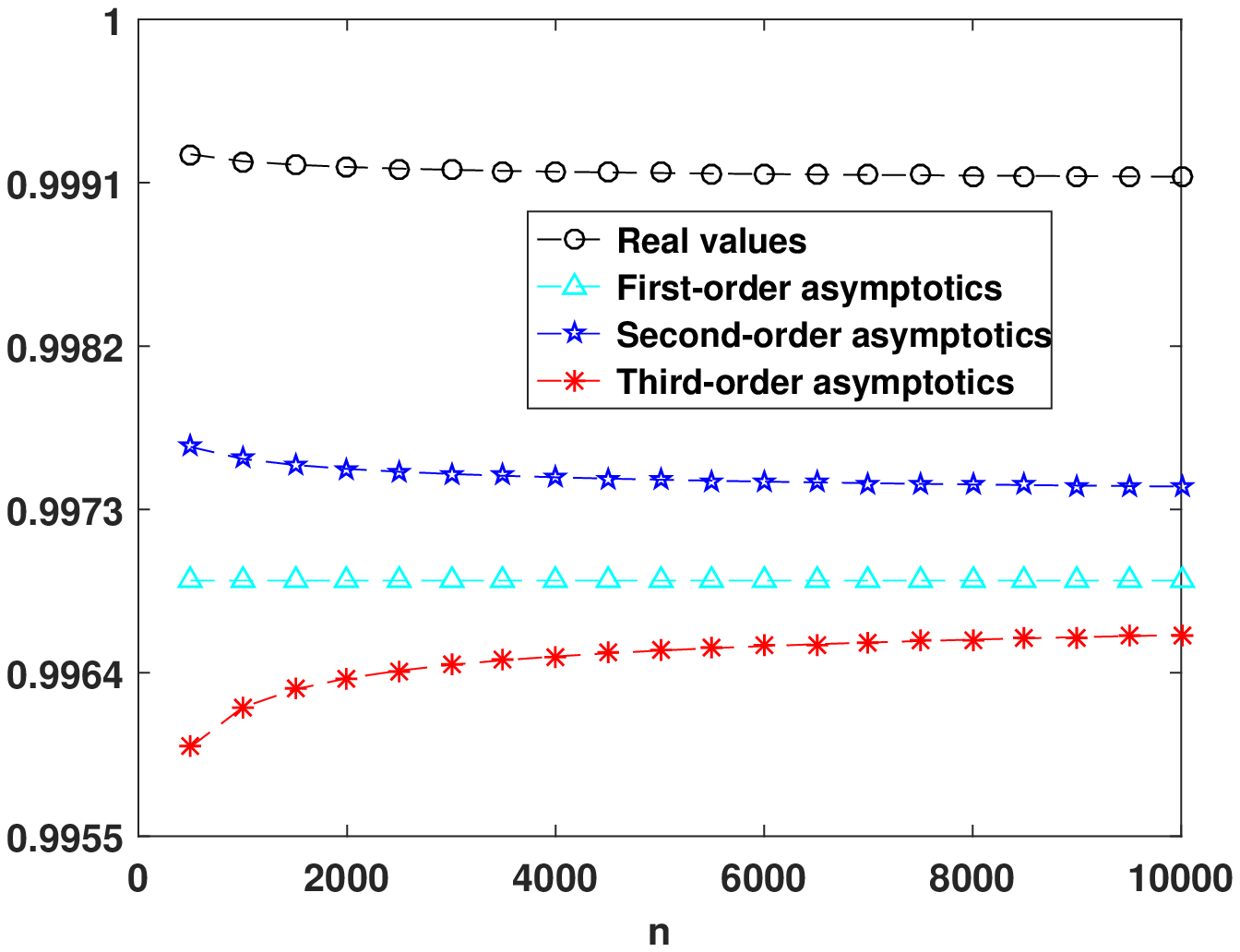}
}
\quad
\subfigure[$(x,y)=(6,2),~\sigma=15$.]{
\includegraphics[scale=0.5]{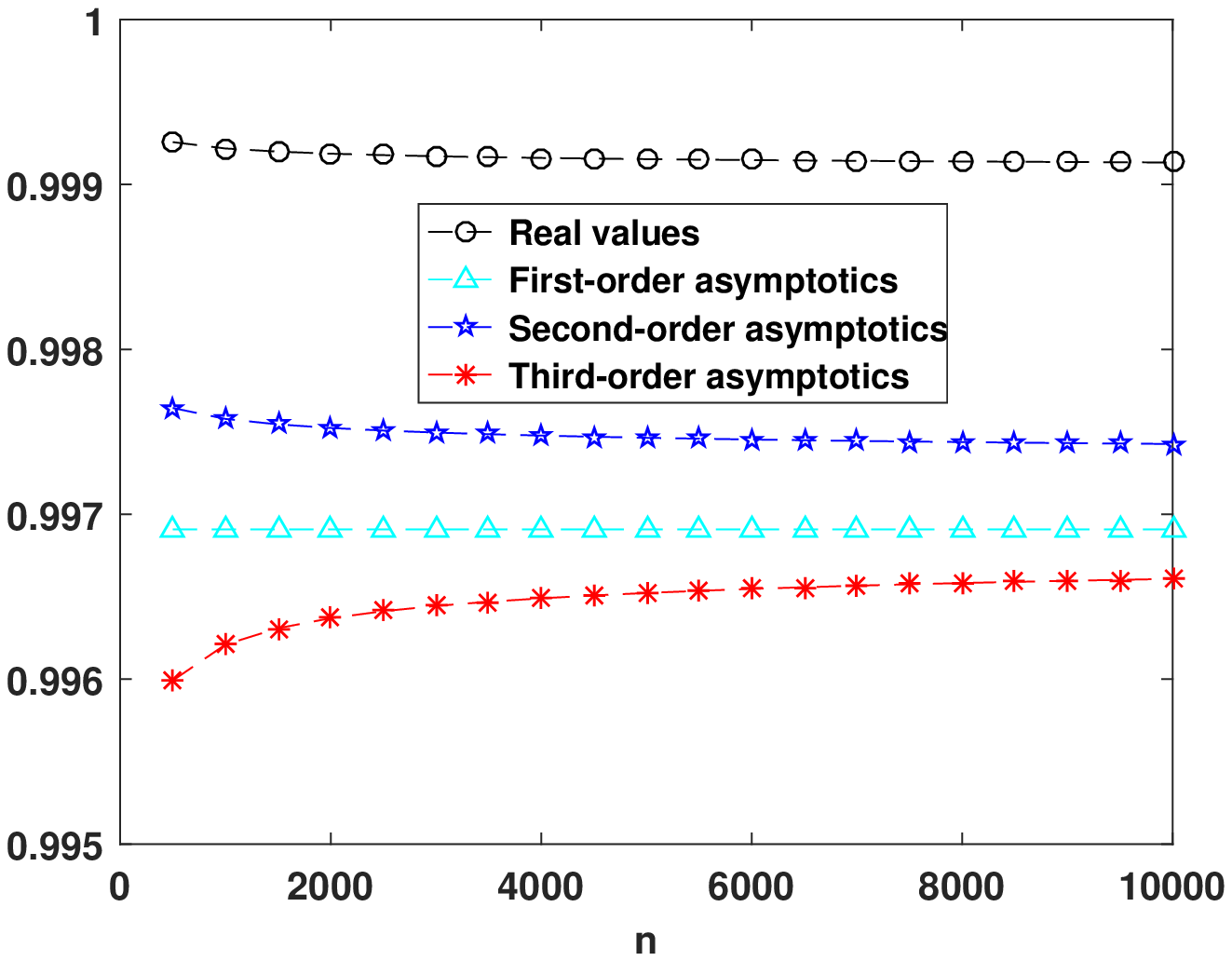}
}
\caption{Real values and its approximations of the joint pdf of $(M_n,~m_n)$ with fixed $(x,y)$ and $k=1$.}\label{fig.1}
\end{figure}

\begin{figure}[htbp]
\centering
\subfigure[$k=0.5$.]{
\includegraphics[scale=0.5]{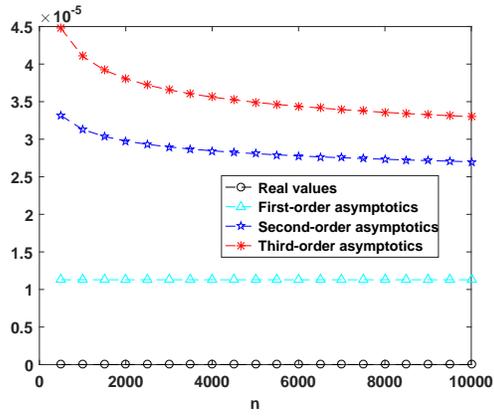}
}
\quad
\subfigure[$k=1.0$.]{
\includegraphics[scale=0.5]{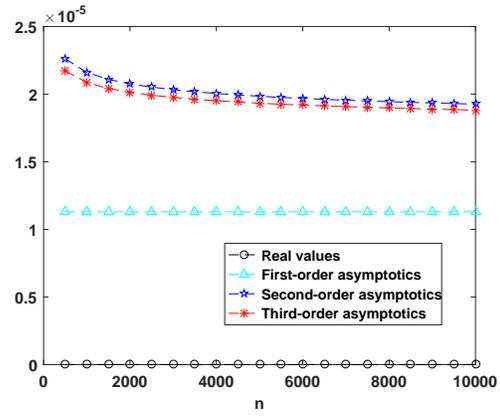}
}
\quad
\subfigure[$k=1.5$.]{
\includegraphics[scale=0.5]{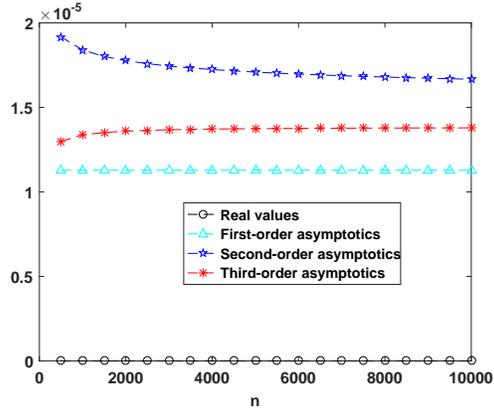}
}
\quad
\subfigure[$k=6.0$.]{
\includegraphics[scale=0.5]{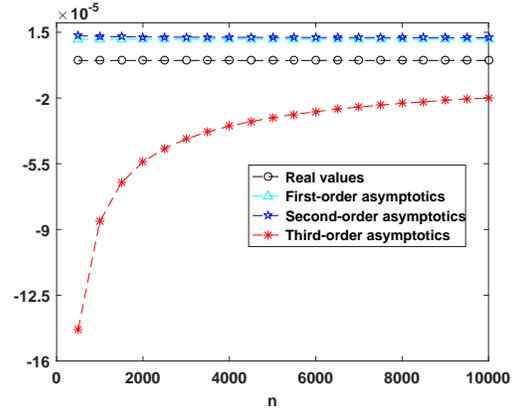}
}
\caption{Real values and its approximations of the joint pdf of $(M_n,~m_n)$ with fixed $(x,y)=(2,6)$ and $\sigma=1$.}\label{fig.2}
\end{figure}

In order to compare the precision of real values with that of its corresponding asymptotics, let
$\Delta_i=|\mathbb{P}(M_n\leq u(x,b_n),~m_n\leq v(y,b_n))-S_i(x,y)|$
and
$\Theta_i=|g_n(x,y)-T_i(x,y)|$, $i=1,2,3$, respectively
stand for the absolute errors. We use Matlable to compute the real values and the asmyptotics of the joint distribution and the joint density of normalized $(M_n,~m_n)$. To demonstrate the precision of all asymptotics with varying $n$, we
then plot the real values and its asymptotics of the joint cdf of the joint pdf with fixed $(x,y)$. Figures \ref{fig.1} and \ref{fig.2} indicate that the larger $n$ the better of all asymptotics.

\begin{figure}[htbp]
\centering
\subfigure[$k=1,~y=2$.]{
\includegraphics[scale=0.25]{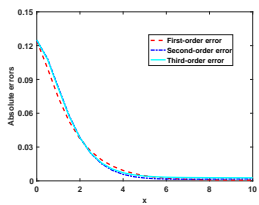}
}
\vspace{0.5cm}
\subfigure[$k=1,~y=10$.]{
\includegraphics[scale=0.25]{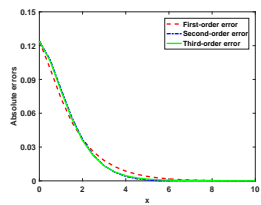}
}
\vspace{0.5cm}
\subfigure[$k=1,~x=2$.]{
\includegraphics[scale=0.25]{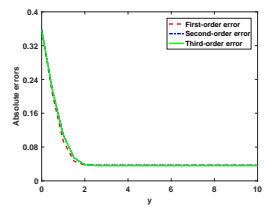}
}
\vspace{0.5cm}
\subfigure[$k=1,~x=10$.]{
\includegraphics[scale=0.25]{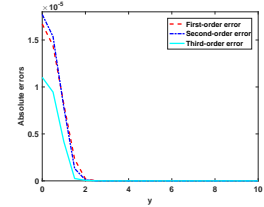}
}
\caption{The absolute errors between real values and its approximations of the joint cdf with $n=500$.}\label{fig.3}
\end{figure}

\begin{figure}[htbp]
\centering
\subfigure[$k=1,~y=2$.]{
\includegraphics[scale=0.25]{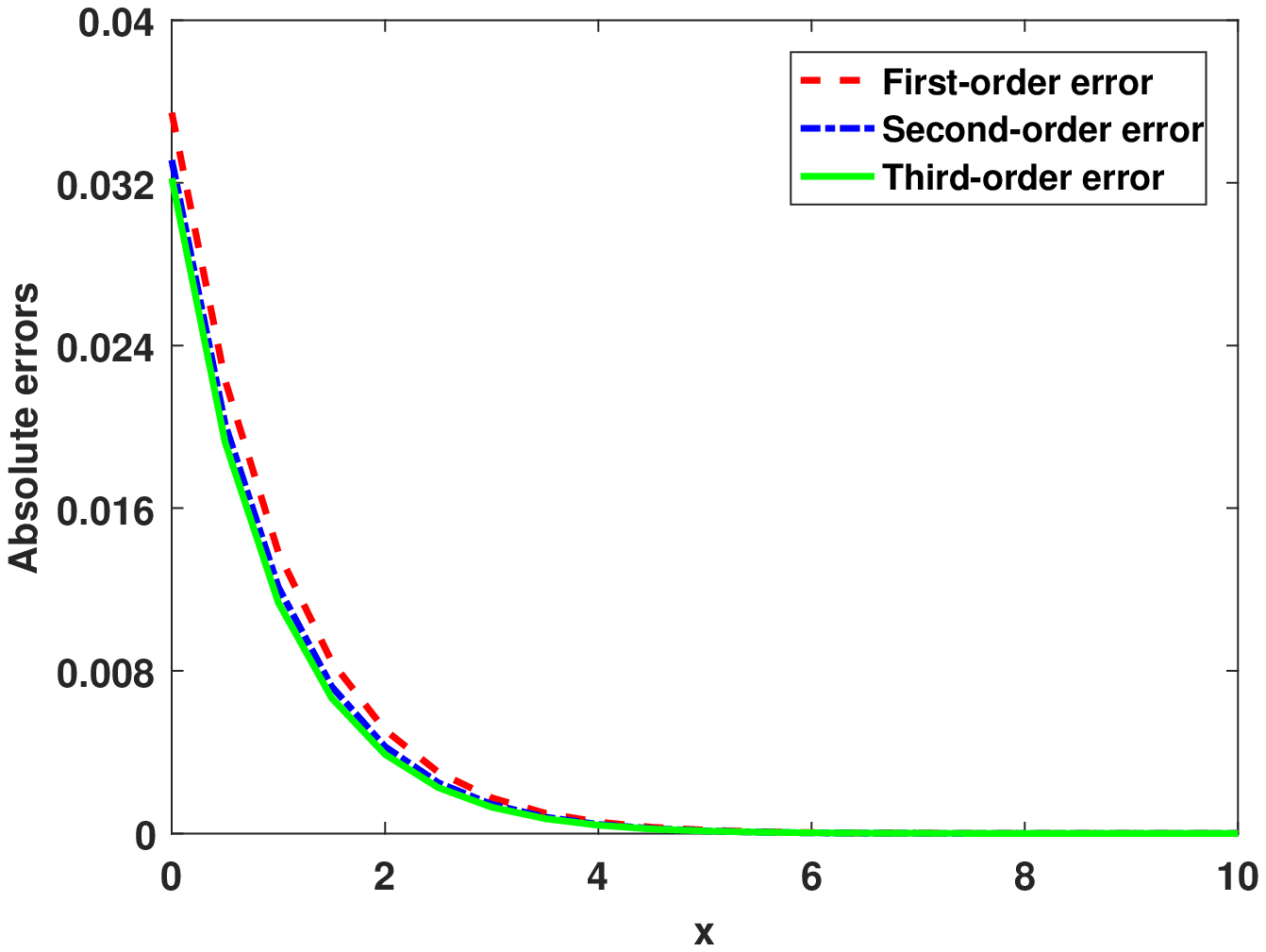}
}
\vspace{0.5cm}
\subfigure[$k=2,~y=2$.]{
\includegraphics[scale=0.25]{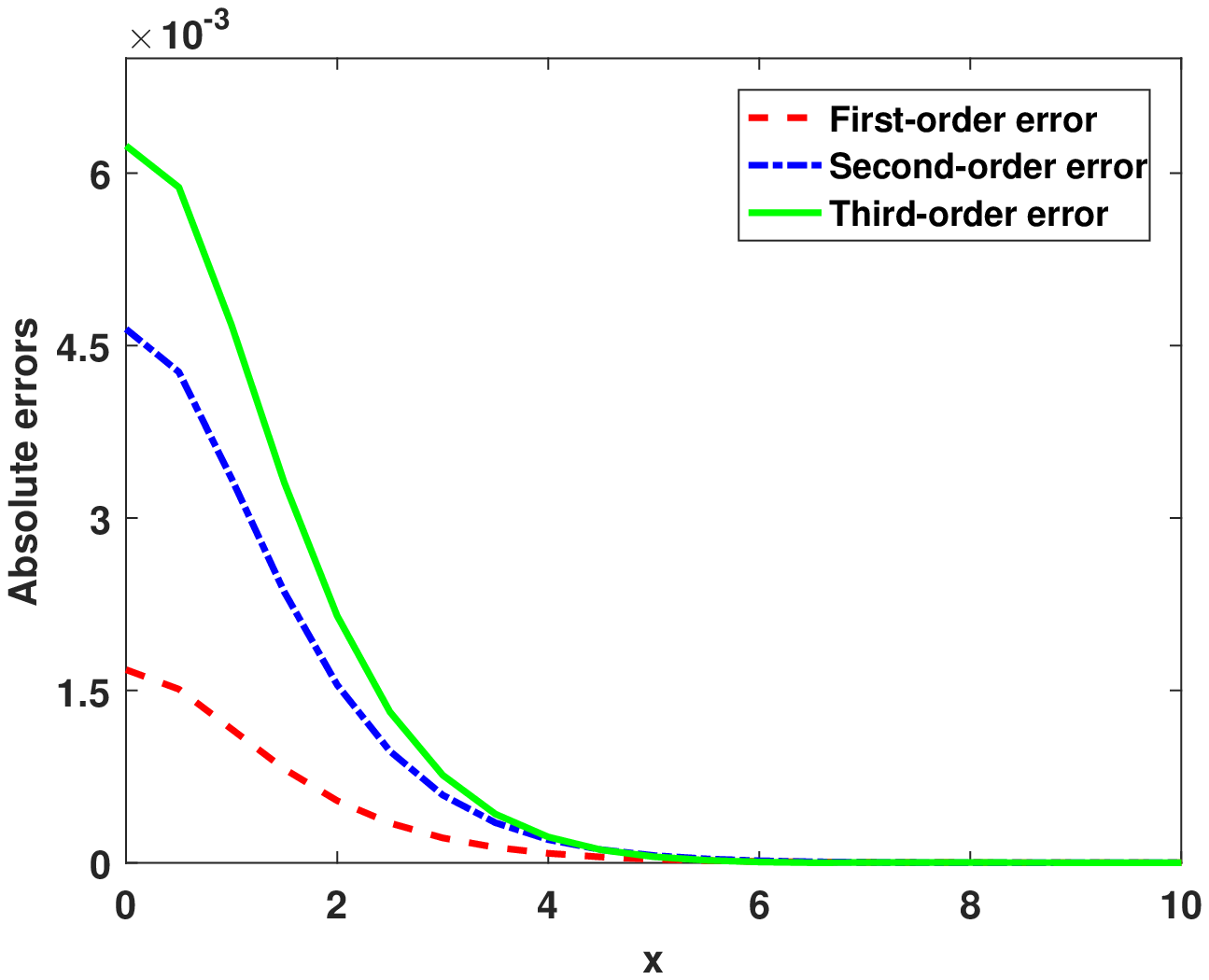}
}
\vspace{0.5cm}
\subfigure[$k=4,~y=2$.]{
\includegraphics[scale=0.25]{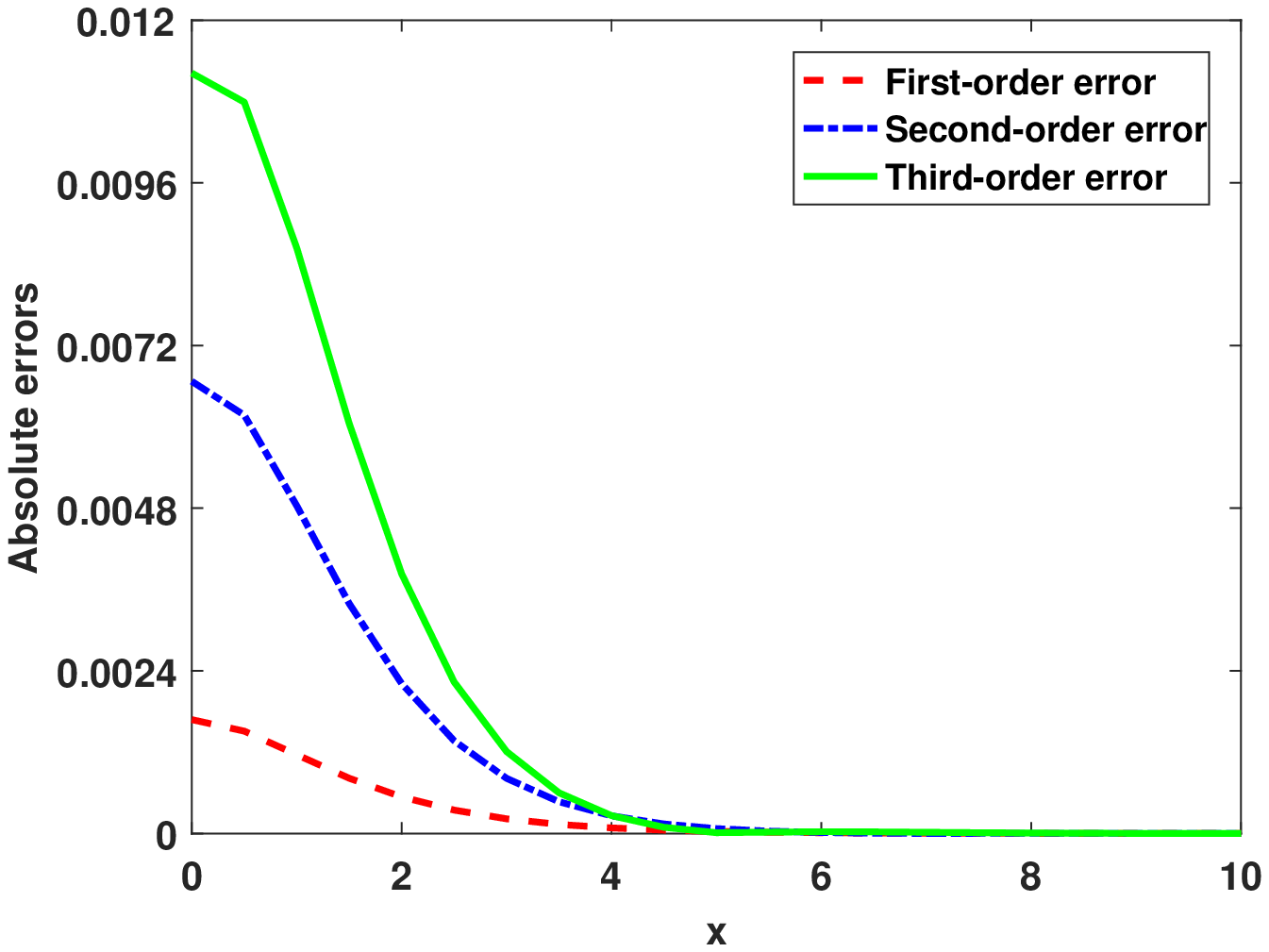}
}
\vspace{0.5cm}
\subfigure[$k=1,~y=10$.]{
\includegraphics[scale=0.25]{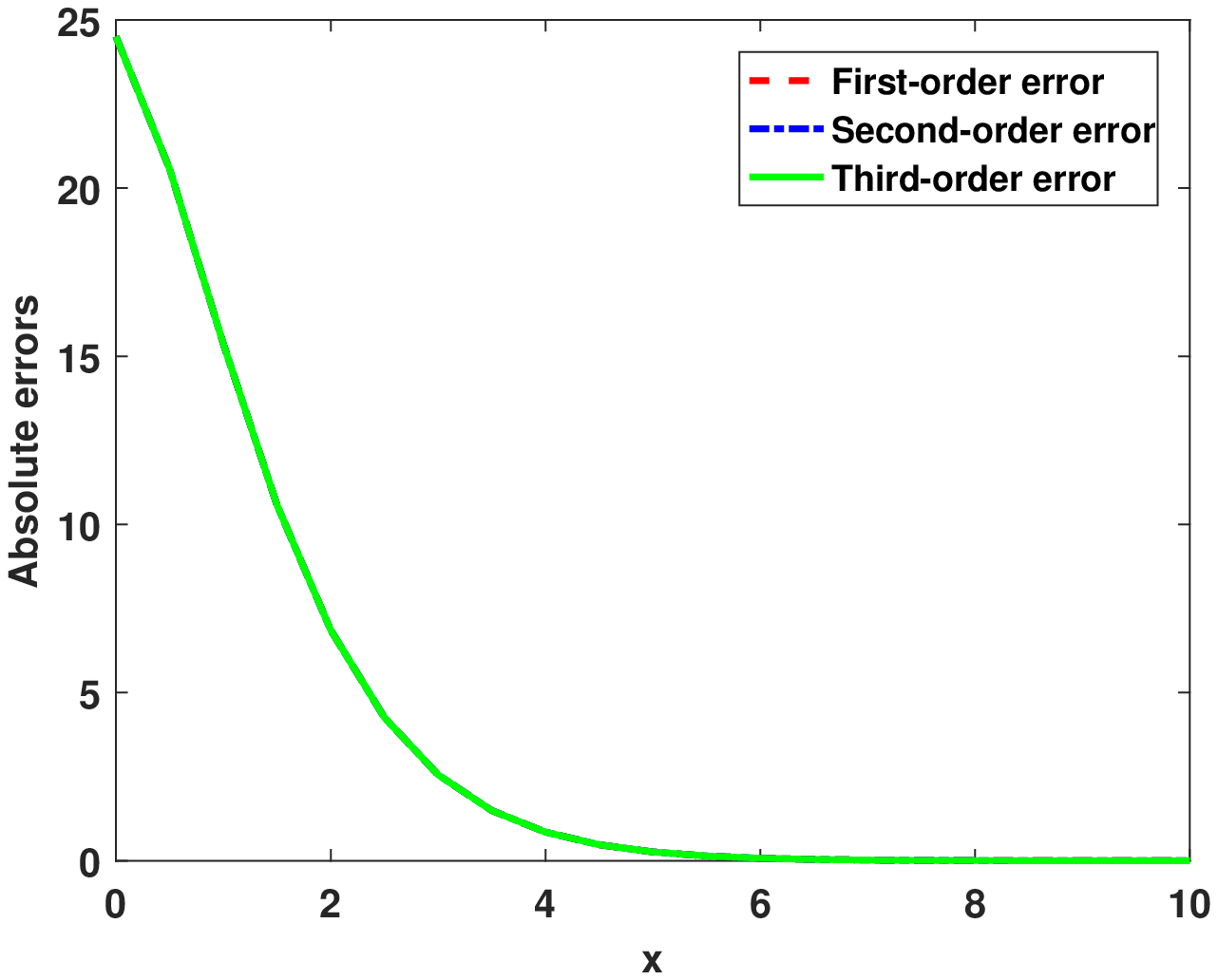}
}
\subfigure[$k=2,~y=10$.]{
\includegraphics[scale=0.25]{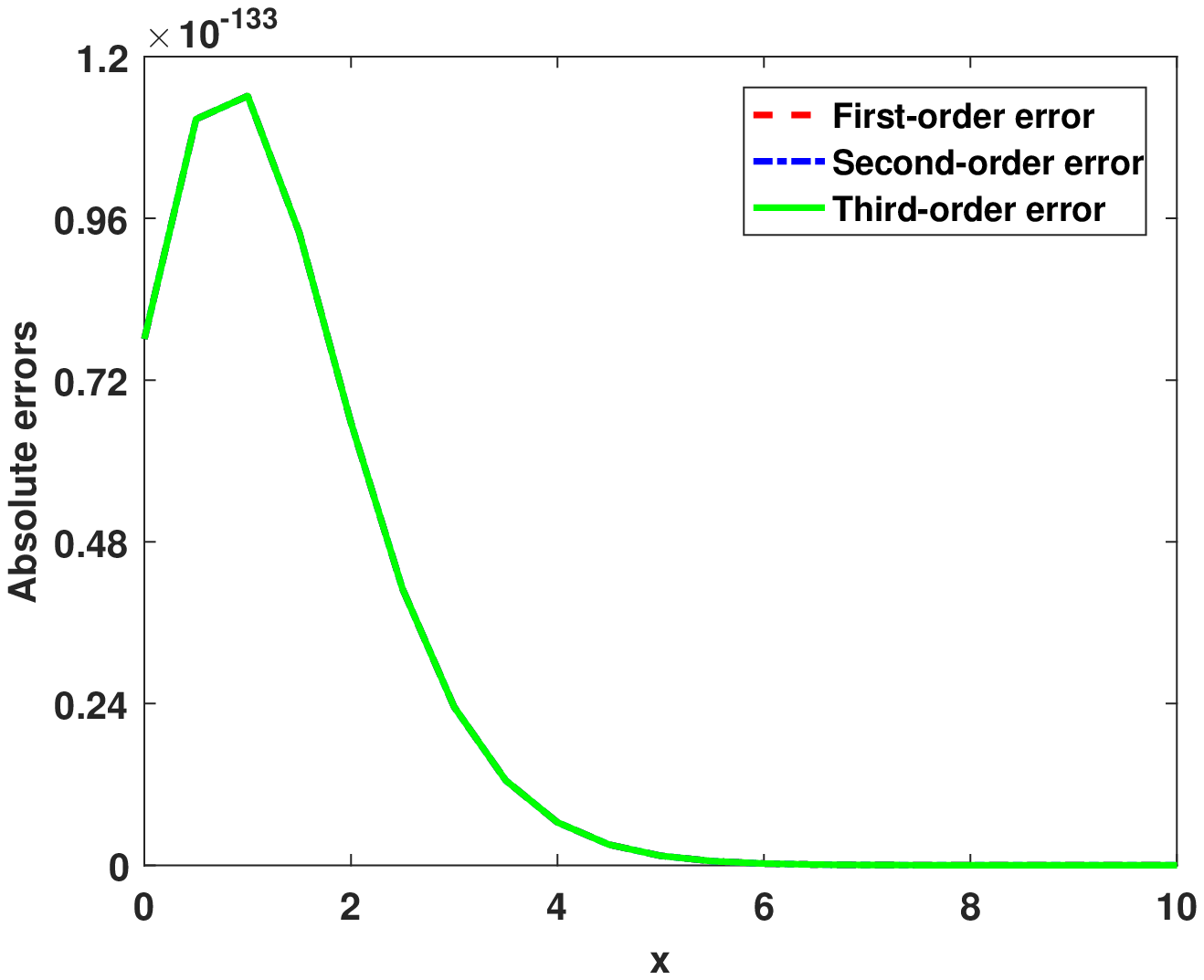}
}
\vspace{0.5cm}
\subfigure[$k=4,~y=10$.]{
\includegraphics[scale=0.25]{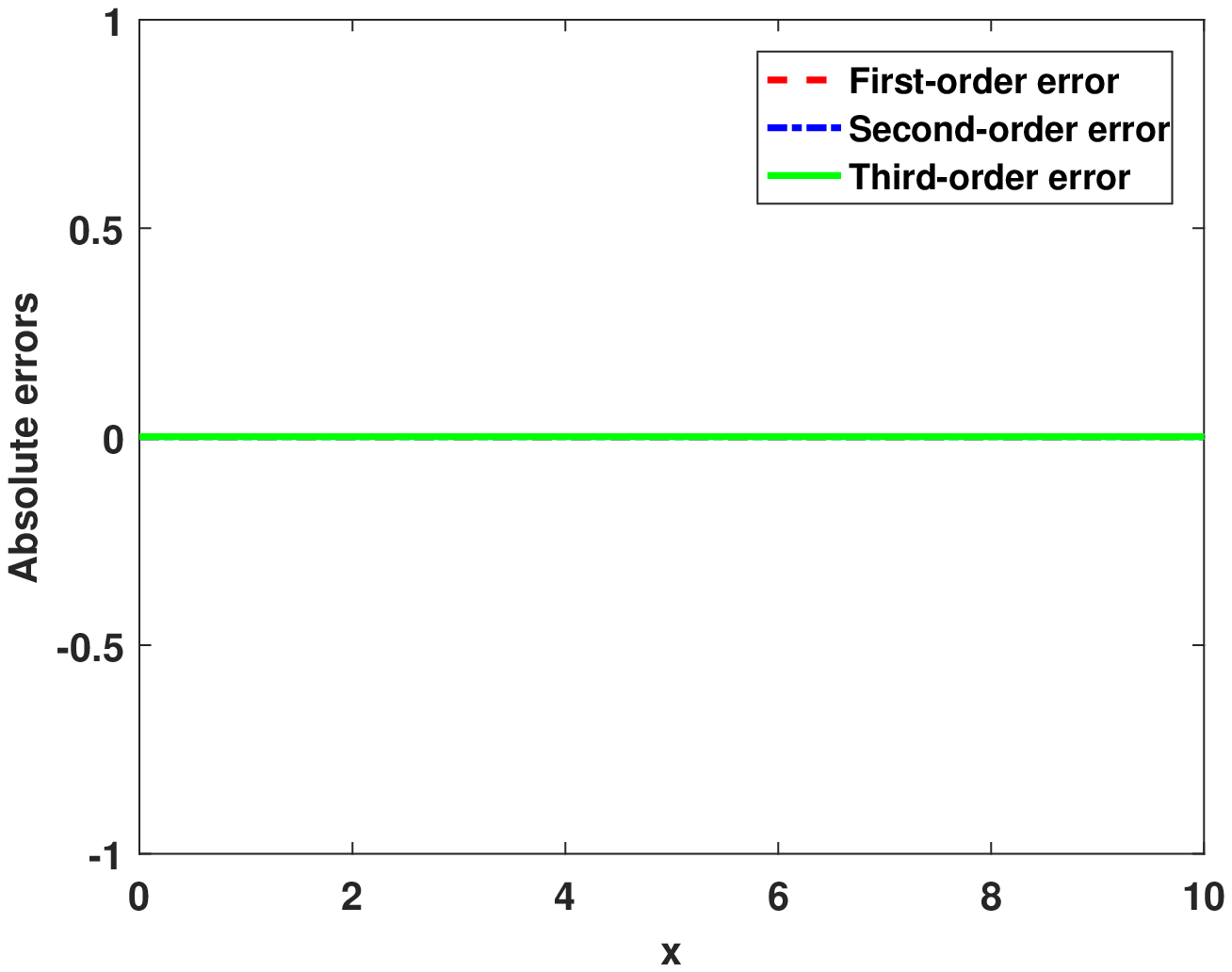}
}
\vspace{0.5cm}
\subfigure[$k=1,~x=2$.]{
\includegraphics[scale=0.25]{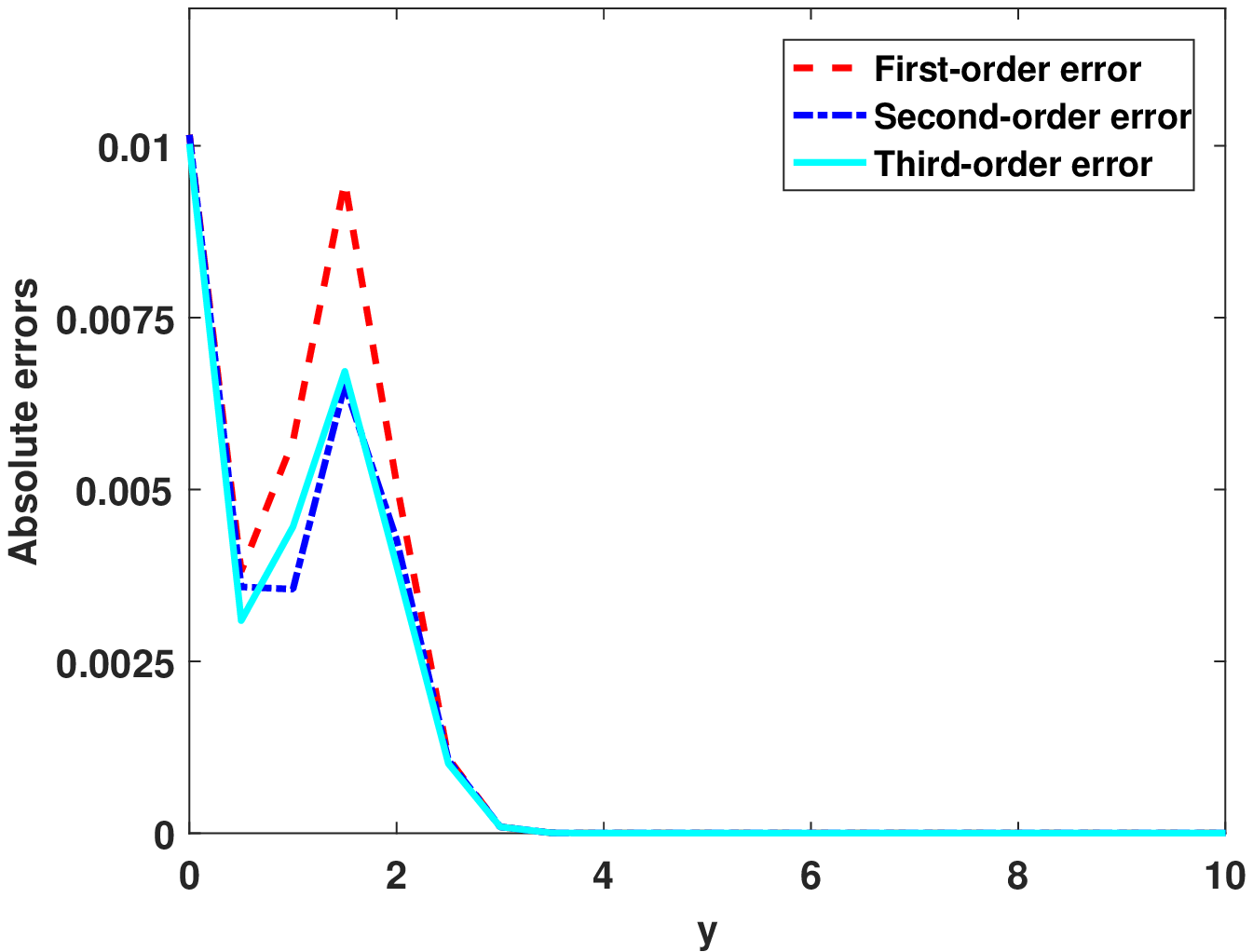}
}
\vspace{0.5cm}
\subfigure[$k=2,~x=2$.]{
\includegraphics[scale=0.25]{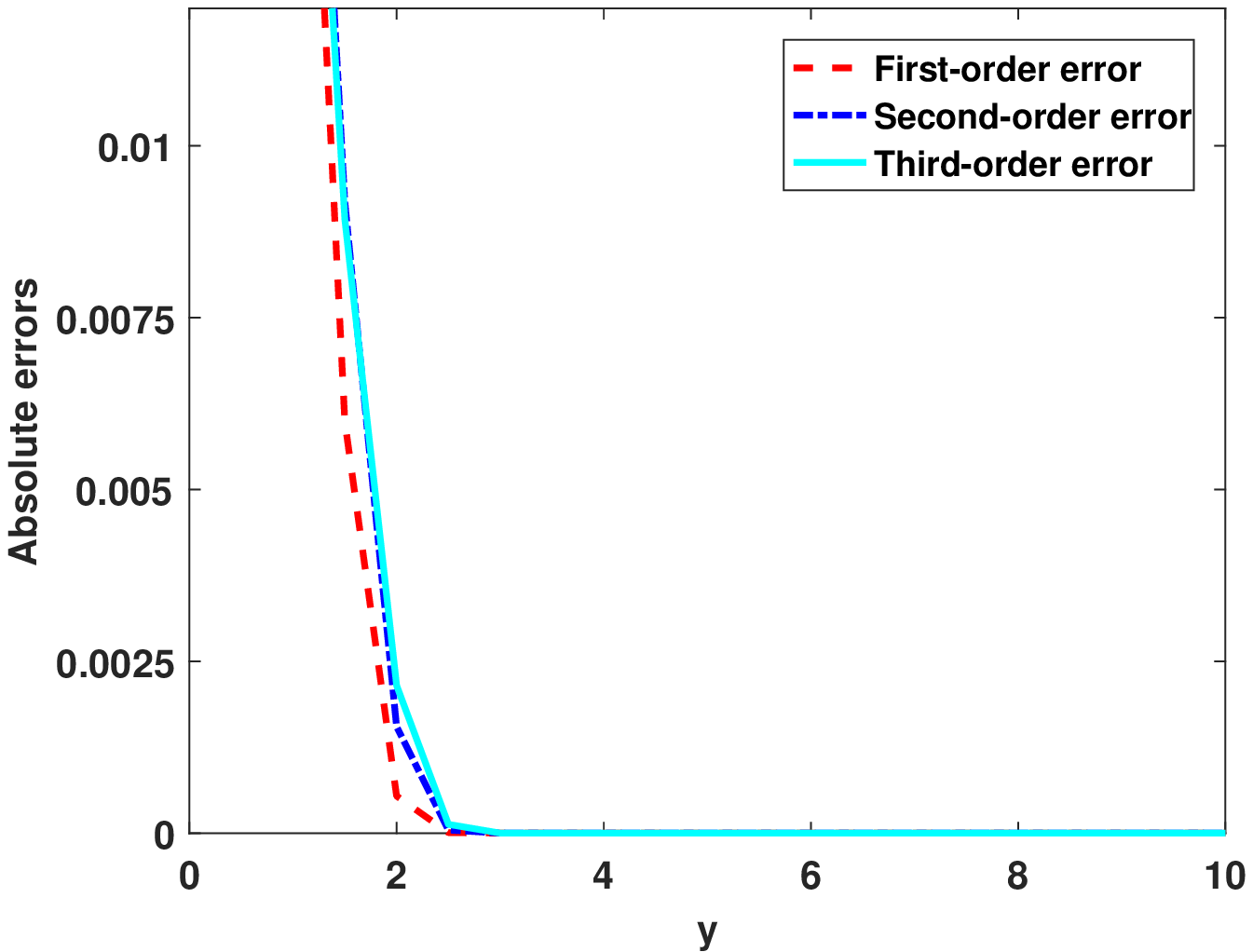}
}
\subfigure[$k=4,~x=2$.]{
\includegraphics[scale=0.25]{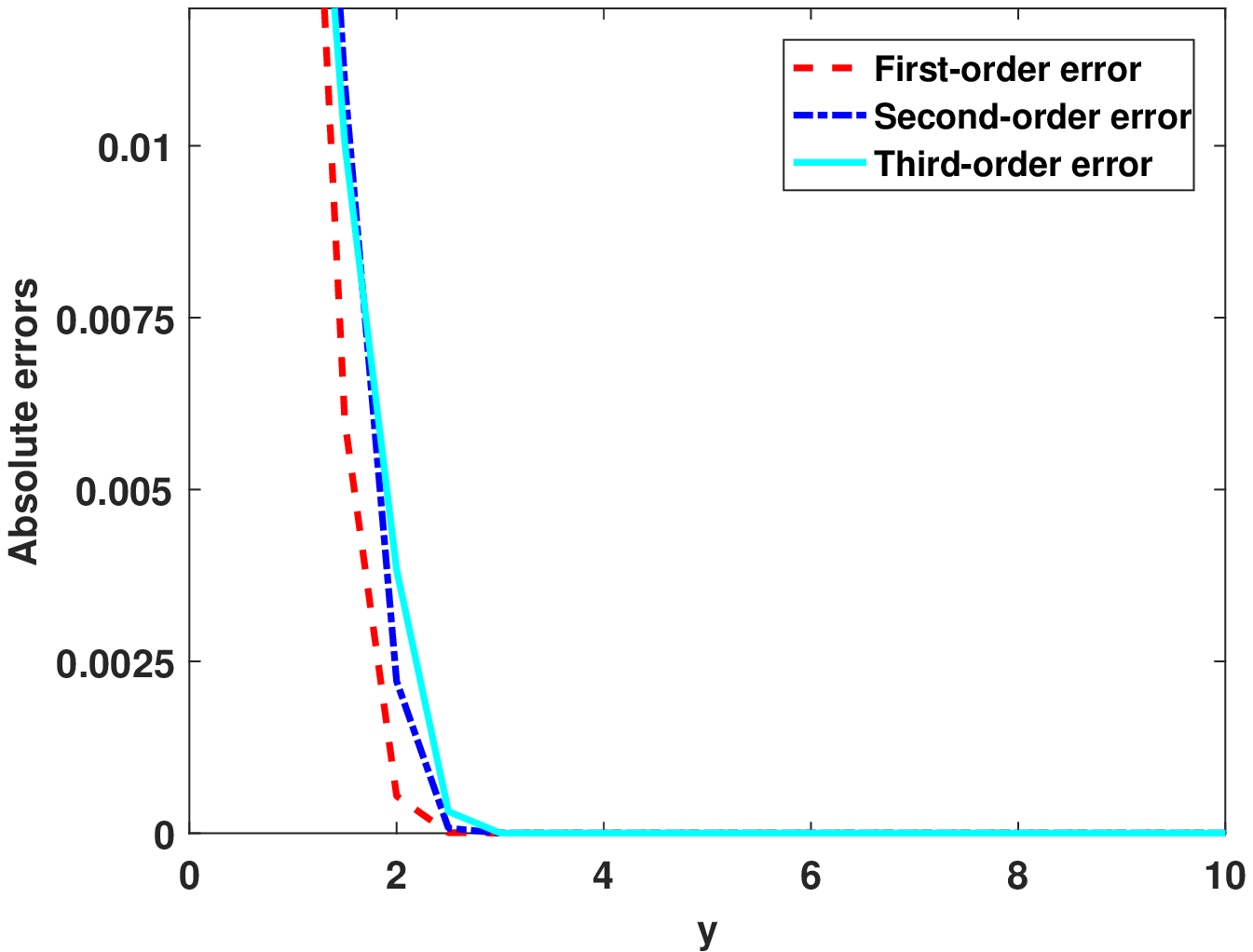}
}
\vspace{0.5cm}
\subfigure[$k=1,~x=10$.]{
\includegraphics[scale=0.25]{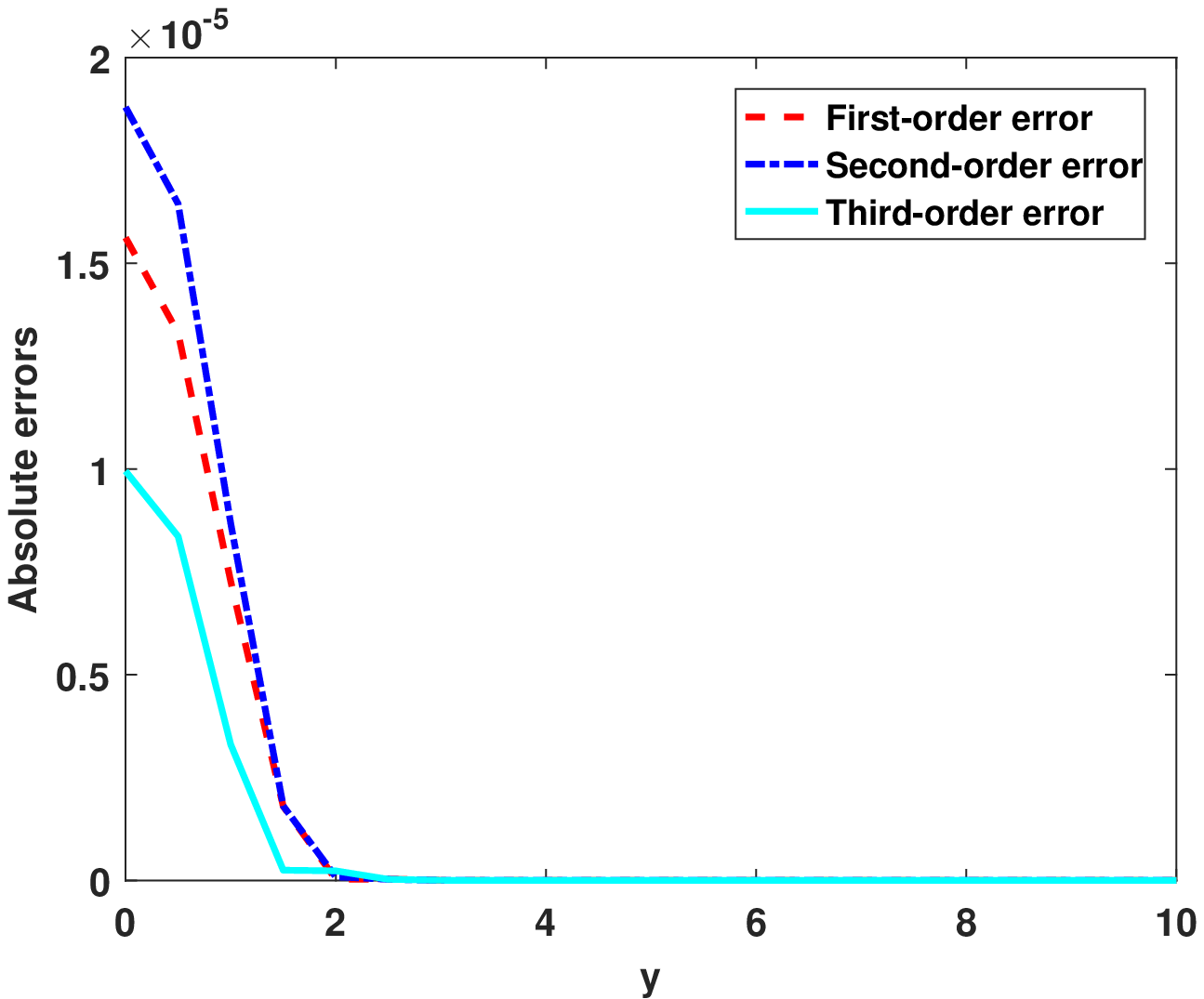}
}
\vspace{0.5cm}
\subfigure[$k=2,~x=10$.]{
\includegraphics[scale=0.25]{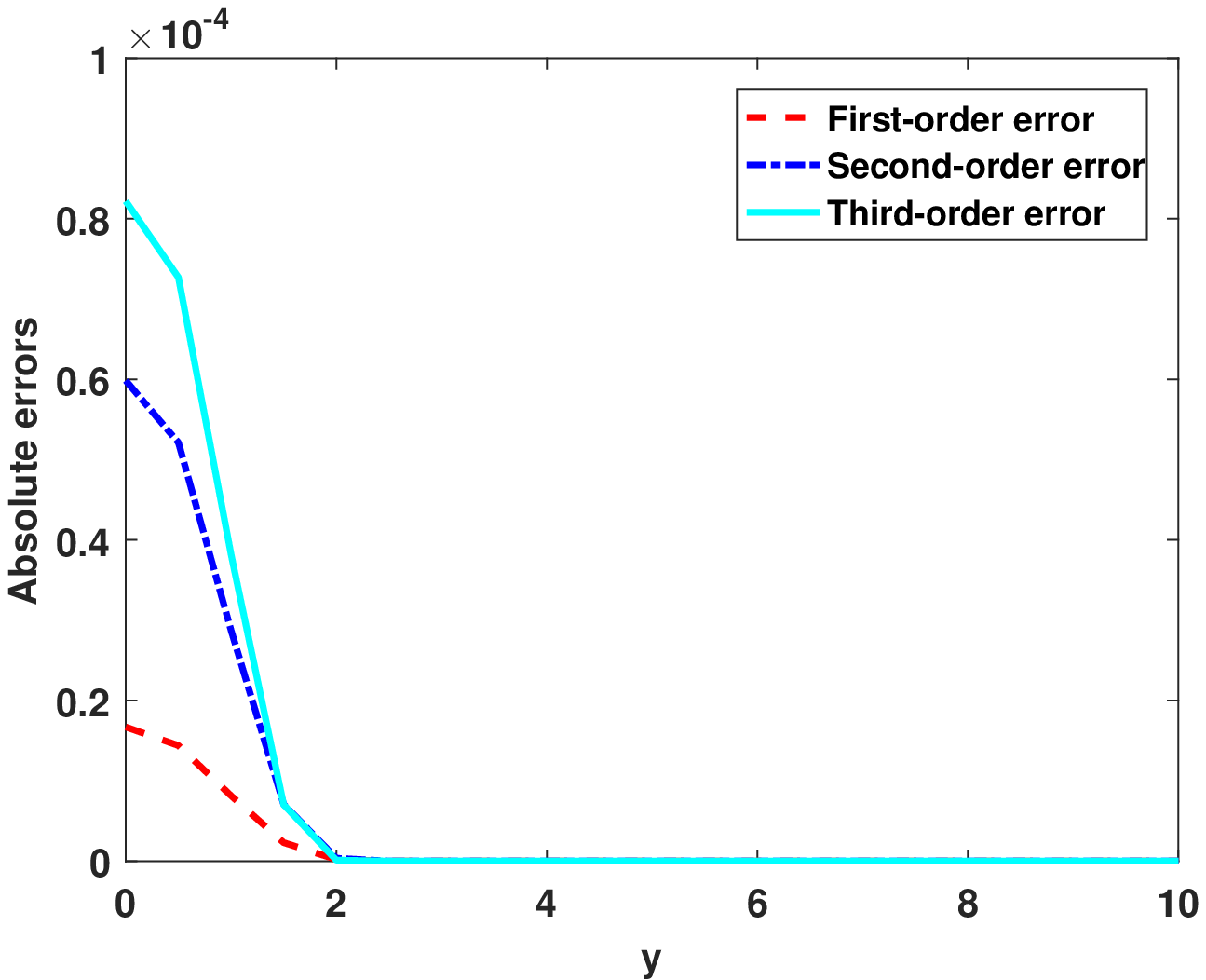}
}
\vspace{0.5cm}
\subfigure[$k=4,~x=10$.]{
\includegraphics[scale=0.25]{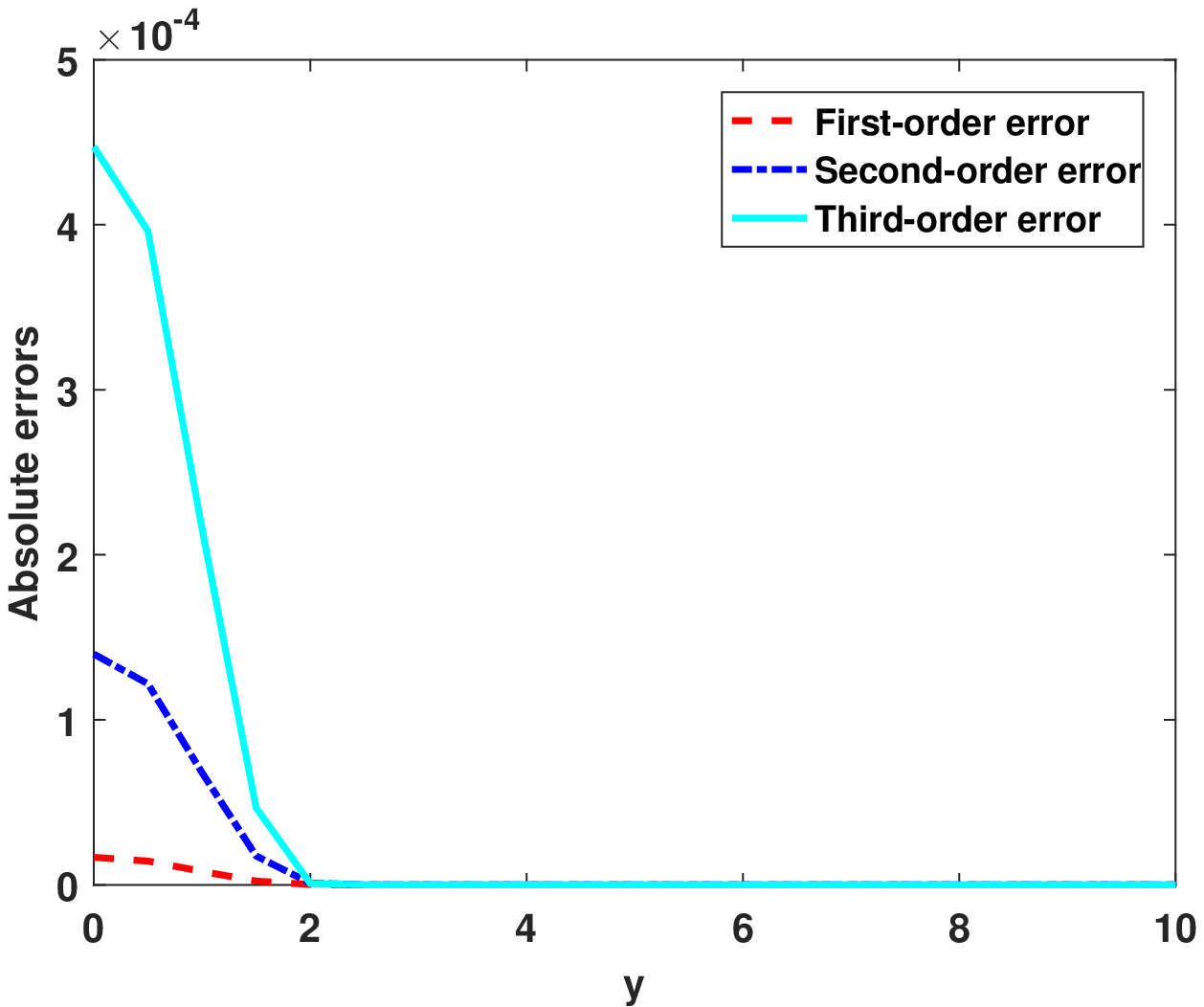}
}
\caption{The absolute errors between real values and its approximations of the joint pdf with $n=500$.}\label{fig.4}
\end{figure}

To evidence the precision of all asymptotics with varying $x$ and $y$, we respectively plot the absolute errors of the joint cdf and the joint pdf with fixed $k$, $\sigma$ and $n=500$ in Figure \ref{fig.3} and Figure \ref{fig.4}.

According to Figures \ref{fig.3} and \ref{fig.4}, we have the following findings: (i) When $x$ increasingly enlarges, all asymptotics of the joint cdf are always closer to its real values except some special cases. (ii) All asymptotics of the joint pdf are more
closer to its real values as $x$ and $y$ becomes large. (iii) The first-order asymptotics may be better except some special cases.

\section{The proofs}
\label{sec.4}
In this section, we present the proofs of main results. In order to prove the main results, we need some necessary lemmas. Lemmas 1 and 2 are given by \cite{Huang and Wang 2018}. The results are expressed as follows.

\begin{lemma}\label{lem.1}
Let $F_k(x)$ and $f_k(x)$ respectively represent the cdf and the pdf of $\mbox{GMD}(k)$, for large $x$, then we get
\begin{align}\label{eq.11}
\notag1-F_k(x)&=f_k(x)\frac{\sigma^2}{k}x^{1-2k}\left[1+\frac{\sigma^2}{k}x^{-2k}
+\frac{1-2k}{k^2}\sigma^4x^{-4k}+O(x^{-6k})\right]\\
&=\frac{\exp(-\frac{1}{2\sigma^2})}{2^{\frac{k}{2}}\sigma^{\frac{1}{k}}\Gamma(1+\frac{k}{2})}
\left[1+\frac{\sigma^2}{k}x^{-2k}
+\frac{1-2k}{k^2}\sigma^4x^{-4k}+O(x^{-6k})\right]\exp\left[-\int^x_1\frac{\tilde{g}(t)}{f(t)}dt\right]
\end{align}
where $f(t)=k^{-1}\sigma^2x^{1-2k}$ and $\tilde{g}(t)=1-k^{-1}\sigma^2x^{-2k}$.
\end{lemma}

\begin{lemma}\label{lem.2}
For $u(x,b_n)$ determined by (\ref{eq.3}) with normalizing constant $b_n$ given by (\ref{eq.4}), we have
\begin{align}\label{eq.12}
b^{2k}_n\left\{b^{2k}_n\left[n\log F_k(u(x,b_n))+e^{-x}\right]-l_k(x)\right\}\to w_k(x),
\end{align}
where $l_k(x)$ and $w_k(x)$ are respectively provided by (\ref{eq.5}) and (\ref{eq.6}).
\end{lemma}

\begin{lemma}\label{lem.3}
Let $h_k(x,y)=n\log[F_k(u(x,b_n))-F_k(v(y,b_n))]+e^{-x}+e^{y}$ with $u(x,b_n)$ and $v(y,b_n)$ determined by (\ref{eq.3}) and (\ref{eq.8}). Then,
\begin{align}\label{eq.13}
b^{2k}_n\left[b^{2k}_nh_k(x,y)-l_k(x,y)\right]\to w_k(x,y),
\end{align}
as $n\to\infty$, where $l_k(x,y)$ and $w_k(x,y)$ are provided by Theorem \ref{the.2}.
\end{lemma}

\noindent
\textbf{Proof.}~~It is easy to verify that
\begin{align}\label{eq.14}
F_k(v(y,b_n)=1-F_k(u(-y,b_n)).
\end{align}
Then,
\begin{align}\label{eq.15}
\notag&\lim_{n\to\infty}b^{2k}_nh_k(x,y)\\
\notag&=\lim_{n\to\infty}nb^{2k}_n\left\{\log\left[1-(1-F_k(u(x,b_n))+F_k(v(y,b_n)))\right]
+n^{-1}(e^{-x}+e^y)\right\}\\
\notag&\overset{\text{(a)}}{=}\lim_{n\to\infty}nb^{2k}_n\bigg\{-\left[1-(1-F_k(u(x,b_n))+F_k(v(y,b_n)))\right]
\\
\notag&\quad -\frac{1}{2}\left[1-(1-F_k(u(x,b_n))+F_k(v(y,b_n)))\right]^2(1+o(1))+n^{-1}(e^{-x}+e^y)\bigg\}\\
\notag&=\lim_{n\to\infty}nb^{2k}_n\bigg\{-\left[1-F_k(u(x,b_n))\right]
-\frac{1}{2}\left[1-F_k(u(x,b_n))\right]^2(1+o(1))+n^{-1}e^{-x}\\
\notag&\quad-F_k(v(y,b_n))
-\frac{1}{2}F_k(v(y,b_n))^2(1+o(1))+n^{-1}e^{y}\\
\notag&\quad-F_k(v(y,b_n))\left[1-F_k(u(x,b_n))\right](1+o(1))\bigg\}\\
\notag&\overset{\text{(b)}}{=}\lim_{n\to\infty}\bigg\{b^{2k}_n[n\log F_k(u(x,b_n))+e^{-x}]+b^{2k}_n[n\log F_k(u(-y,b_n))+e^{y}]+O(n^{-1}b^{2k}_n)\bigg\}\\
\notag&\overset{\text{(c)}}{=}l_k(x)+l_k(-y)\\
&:=l_k(x,y),
\end{align}
where (a) and (b) follow from $\log(1-x)=-x-x^2/2-x^3/3+o(x^3),~x\to 0$, and (c) is due to Lemma \ref{lem.2} and (\ref{eq.14}).

Likewise, by Lemma \ref{lem.2}, we can obtain
\begin{align}\label{eq.16}
\notag&\lim_{n\to\infty}b^{2k}_n\left[b^{2k}_nh_k(x,y)-l_k(x,y)\right]\\
\notag&=\lim_{n\to\infty}\bigg\{b^{2k}_n\left\{b^{2k}_n\left[n\log F_k(u(x,b_n))+e^{-x}\right]-l_k(x)\right\}\\
\notag&\quad +b^{2k}_n\left\{b^{2k}_n\left[n\log F_k(u(-y,b_n))+e^{y}\right]-l_k(-y)\right\}+O(n^{-1}b^{-4k}_n)\bigg\}\\
\notag&=w_k(x)+w_k(-y)\\
&:=w_k(x,y).
\end{align}
The proof is complete.

\qed

\noindent
\textbf{Proof of Proposition \ref{pro.1}.}~~It is not hard to check that $m_n=-\max_{1\leq i\leq n}\{-X_i\}$ with $-X_i$ following $\mbox{GMD}(k)$. Combining with (\ref{eq.2}) and (\ref{eq.3}), it implies that
\begin{align}\label{eq.17}
\lim_{n\to\infty}\mathbb{P}(m_n\leq v(y,b_n))=1-\Lambda(-y),
\end{align}
where
\begin{align}\notag
v(y,b_n)=k^{-1}\sigma^2b^{1-2k}_ny-b_n.
\end{align}
By (\ref{eq.14}), Lemmas \ref{lem.1} and \ref{lem.2} and analogous discussion for the proof of Theorem $2.1$ \cite{Huang and Wang 2018}, the desired result can be derived.

\qed

\noindent
\textbf{Proof of Theorem \ref{the.1}.}~~By (\ref{eq.2}), (\ref{eq.17}) together with Theorem $1.8.2$ \cite{Leadbetter 1983}, the claim is deduced.

\qed

\noindent
\textbf{Proof of Theorem \ref{the.2}.}~~By (\ref{eq.15}), it follows that
\begin{align}\label{eq.18}
h_k(x,y)\to 0~\mbox{and}~\left|\sum^{\infty}_{j=3}\frac{h^{j-3}_k(x,y)}{j!}\right|\leq\exp(|h_k(x,y)|)\to 1,
\end{align}
as $n\to\infty$. By (\ref{eq.15}), (\ref{eq.16}) and (\ref{eq.18}), we get
\begin{align}\label{eq.19}
\notag &\lim_{n\to\infty}b^{2k}_n\left\{b^{2k}_n\left[\mathbb{P}(M_n\leq u(x,b_n),~m_n>v(y,b_n))-\Lambda(x)\Lambda(-y)\right]-l_k(x,y)\Lambda(x)\Lambda(-y)\right\}\\
\notag&=\lim_{n\to\infty}b^{2k}_n\left\{b^{2k}_n\left[\exp(h_k(x,y))-1\right]-l_k(x,y)\right\}\Lambda(x)\Lambda(-y)\\
\notag&\overset{\text{(a)}}{=}\lim_{n\to\infty}\left\{b^{2k}_n[b^{2k}_nh_k(x,y)-l_k(x,y)]+b^{4k}_nh^2_k(x,y)\left[\frac{1}{2}
+h_k(x,y)\sum^{\infty}_{j=3}\frac{h^{j-3}_k(x,y)}{j!}\right]\right\}\Lambda(x)\Lambda(-y)\\
&=\left[w_k(x,y)+\frac{1}{2}l^2_k(x,y)\right]\Lambda(x)\Lambda(-y),
\end{align}
where (a) is because of $e^x=1+x+x^2(1+O(x)),~x\to 0$. By (\ref{eq.28}) and (\ref{eq.19}), we get
\begin{align}\label{eq.20}
\notag&\mathbb{P}(M_n\leq u(x,b_n),~m_n\leq v(y,b_n))\\
\notag&=\mathbb{P}(M_n\leq u(x,b_n))-\mathbb{P}(M_n\leq u(x,b_n),~m_n>v(y,b_n))\\
\notag&=\Lambda(x)(1-\Lambda(-y))+b^{-2k}_n[l_k(x)\Lambda(x)-l_k(x,y)\Lambda(x)\Lambda(-y)]\\
&\quad+b^{-4k}_n\left\{\left[w_k(x)+\frac{1}{2}l_k^2(x)\right]\Lambda(x)
-\left[w_k(x,y)+\frac{1}{2}l_k^2(x,y)\right]\Lambda(x)\Lambda(-y)\right\}+O(b^{-6k}_n),
\end{align}
as $n$ is large. We complete the proof.

\qed

\noindent
\textbf{Proof of Theorem \ref{the.3}.}~~It is known that $n(1-F_k(u(x,b_n)))\to e^{-x}$ for large $n$ implies $1-F_k(u(x,b_n))=O(n^{-1})$, and similarly, $F_k(v(y,b_n))=O(n^{-1})$. Therefore,
\begin{align}\label{eq.21}
\notag&[F_k(u(x,b_n))-F_k(v(y,b_n))]^{n-2}\\
\notag&=\frac{[F_k(u(x,b_n))-F_k(v(y,b_n))]^n}{\{1-\{1-[F_k(u(x,b_n))-F_k(v(y,b_n))]\}\}^2}\\
\notag&\overset{\text{(a)}}{=}(1+O(n^{-1}))[F_k(u(x,b_n))-F_k(v(y,b_n))]^n\\
&\overset{\text{(b)}}{=}\Lambda(x)\Lambda(-y)+b^{-2k}_nl_k(x,y)\Lambda(x)\Lambda(-y)
+b^{-4k}_n\left[w_k(x,y)+\frac{1}{2}l_k^2(x,y)\right]\Lambda(x)\Lambda(-y)+O(b^{-6k}_n)
\end{align}
for large $n$, where (a) is from $1/(1-x)=1+x+x^2+x^3(1+o(1)),~x\to 0$ and (b) is thanks to (\ref{eq.19}).

Employing the following Taylor's expansion
\begin{align}
\notag(1+x)^{\alpha}=1+\alpha x+\frac{\alpha(\alpha-1)}{2}x^2
+\frac{\alpha(\alpha-1)(\alpha-2)}{6}x^3+O(x^4),~x\to 0,~\alpha\in\mathbb{R},
\end{align}
and the fact in view of (\ref{eq.4}) that $b^{2k}_n\sim\log n$ for large $n$, we get
\begin{align}\label{eq.22}
\notag&\left(\frac{\sigma^2}{k}b^{1-2k}_nx+b_n\right)^{2k}\\
\notag&=b^{2k}_n\left(1+\frac{\sigma^2}{k}b^{-2k}_nx\right)^{2k}\\
&=b^{2k}_n\left[1+2\sigma^2b^{-2k}_nx+\frac{2k-1}{k}\sigma^4b^{-4k}_nx^2
+\frac{(2k-1)(2k-2)}{3k^2}\sigma^6b^{-6k}_nx^3+O(b^{-8k}_n)\right].
\end{align}
Applying (\ref{eq.22}), we get
\begin{align}\label{eq.23}
\notag&\exp\left(-\frac{1}{2\sigma^2}\left(\frac{\sigma^2}{k}b^{1-2k}_nx+b_n\right)^{2k}
+\frac{1}{2\sigma^2}b^{2k}_n+x\right)\\
\notag&=\exp\left(-\frac{2k-1}{2k}\sigma^2b^{-2k}_nx^2
-\frac{(2k-1)(k-1)}{3k^2}\sigma^4b^{-4k}_nx^3+O(b^{-6k}_n)\right)\\
&\overset{\text{(a)}}{=}1-\frac{2k-1}{2k}\sigma^2b^{-2k}_nx^2+\frac{(2k-1)\sigma^4}{k^2}\left(\frac{2k-1}{8}x^4
-\frac{k-1}{3}x^3\right)b^{-4k}_n+O(b^{-6k}_n),
\end{align}
where (a) follows from the fact that $e^x=1+x+x^2/2+x^3/6+O(x^4),~x\to 0$.

Utilizing the following Taylor's expansion
\begin{align}
\notag(1-x)^{-1}=1+x+x^2+x^3+O(x^4),~x\to 0,
\end{align}
we get
\begin{align}\label{eq.24}
\notag&\left[1+\frac{\sigma^2}{k}b^{-2k}_n+\frac{1-2k}{k^2}\sigma^4b^{-4k}_n+O(b^{-6k}_n)\right]^{-1}\\
&=1-\frac{\sigma^2}{k}b^{-2k}_n+\frac{2\sigma^4}{k}b^{-4k}_n+O(b^{-6k}_n).
\end{align}
Hence, by (\ref{eq.4}), (\ref{eq.22}), (\ref{eq.23}), (\ref{eq.24}) and Lemma \ref{lem.1}, we get
\begin{align}\label{eq.25}
\notag&\frac{\sigma^2}{k}\frac{nf_k(u(x,b_n))}{b^{2k-1}_n}e^x\\
\notag&=\frac{\sigma^2}{k}\frac{k}{2^{k/2}\sigma^{2+1/k}\Gamma(1+k/2)}
\frac{(\frac{\sigma^2}{k}b^{1-2k}_nx+b_n)^{2k}
\exp(-\frac{1}{2\sigma^2}(\frac{\sigma^2}{k}b^{1-2k}_nx+b_n)^{2k})}{b^{2k-1}_n(1-F_k(b_n))}e^x\\
\notag&=\frac{(\frac{\sigma^2}{k}b^{1-2k}_nx+b_n)^{2k}
\exp(-\frac{1}{2\sigma^2}(\frac{\sigma^2}{k}b^{1-2k}_nx+b_n)^{2k}+x)}{b^{2k}_n\exp(-\frac{1}{2\sigma^2}b^{2k}_n)}
[1+\frac{\sigma^2}{k}b^{-2k}_n+\frac{1-2k}{k^2}\sigma^4b^{-4k}_n+O(b^{-6k}_n)]^{-1}\\
\notag&=1-(\frac{2k-1}{2k}x^2-2x+\frac{1}{k})\sigma^2b^{-2k}_n+\frac{\sigma^4}{k}
[\frac{(2k-1)^2}{8k}x^4-\frac{(2k+1)(4k-1)}{3k}x^3\\
&\quad+\frac{(2k+1)(2k-1)}{2k}x^2-2x+2]b^{-4k}_n+O(b^{-6k}_n)
\end{align}
and
\begin{align}\label{eq.26}
\notag&\frac{\sigma^2}{k}\frac{nf_k(v(y,b_n))}{b^{2k-1}_n}e^{-y}\\
\notag&=\frac{\sigma^2}{k}\frac{k}{2^{k/2}\sigma^{2+1/k}\Gamma(1+k/2)}
\frac{(\frac{\sigma^2}{k}b^{1-2k}_n(-y)+b_n)^{2k}
\exp(-\frac{1}{2\sigma^2}(\frac{\sigma^2}{k}b^{1-2k}_n(-y)+b_n)^{2k})}{b^{2k-1}_n(1-F_k(b_n))}e^{-y}\\
\notag&=1-(\frac{2k-1}{2k}y^2+2y+\frac{1}{k})\sigma^2b^{-2k}_n+\frac{\sigma^4}{k}
[\frac{(2k-1)^2}{8k}y^4+\frac{(2k+1)(4k-1)}{3k}y^3\\
&\quad+\frac{(2k+1)(2k-1)}{2k}y^2+2y+2]b^{-4k}_n+O(b^{-6k}_n)
\end{align}
for large $n$. Combining with (\ref{eq.21}), (\ref{eq.25}) and (\ref{eq.26}), we have
\begin{align}
\notag&\Delta_n(g_n,g;x,y)\\
\notag&=\Lambda(x)\Lambda(-y)e^{-x}e^y\bigg\{(1-n^{-1})
\frac{\frac{\sigma^2}{k}nf_k(u(x,b_n))e^x}{b^{2k-1}_n}\frac{\frac{\sigma^2}{k}nf_k(v(y,b_n))e^{-y}}{b^{2k-1}_n}\\
\notag&\quad\times\frac{[F_k(u(x,b_n))-F_k(v(y,b_n))]^{n-2}}{\Lambda(x)\Lambda(-y)}-1\bigg\}\\
\notag&=\Lambda(x)\Lambda(-y)e^{-x}e^y\bigg\{\left[l_k(x,y)
-\sigma^2\left[\frac{2k-1}{2k}(x^2+y^2)-2(x-y)+\frac{2}{k}\right]\right]b^{-2k}_n\\
\notag&\quad+\bigg\{\sigma^4\bigg\{\frac{1}{k}\bigg[\frac{(2k-1)^2}{8k}(x^4+y^4)-\frac{(2k+1)(4k-1)}{3k}(x^3-y^3)
+\frac{(2k+1)(2k-1)}{2k}(x^2+y^2)\\
\notag&\quad-2(x-y)+4\bigg]+\left(\frac{2k-1}{2k}x^2-2x+\frac{1}{k}\right)\left(\frac{2k-1}{2k}x^2-2x+\frac{1}{k}\right)
\bigg\}\\
\notag&\quad-\sigma^2\left[\frac{2k-1}{2k}(x^2+y^2)-2(x-y)+\frac{2}{k}\right]l_k(x,y)+w_k(x,y)
+\frac{1}{2}l^2_k(x,y)\bigg\}b^{-4k}_n+O(b^{-6k}_n)\bigg\}\\
\notag&:=\Lambda(x)\Lambda(-y)e^{-x}e^y\left[C_1(x,y)b^{-2k}_n+C_2(x,y)b^{-4k}_n+O(b^{-6k}_n)\right]
\end{align}
for large $n$. We complete the proof.

\qed

\noindent
{\bf Acknowledgements}~~

\end{document}